\documentclass[twoside,english]{elsarticle}
\usepackage[T1]{fontenc}
\usepackage[utf8]{inputenc}
\usepackage{amsmath}
\usepackage{amsthm}
\usepackage{amssymb}

\makeatletter
\theoremstyle{plain}
\newtheorem{thm}{\protect\theoremname}[section]
  \theoremstyle{plain}
  \newtheorem{assumption}[thm]{\protect\assumptionname}
  \theoremstyle{plain}
  \newtheorem{cor}[thm]{\protect\corollaryname}
  \theoremstyle{definition}
  \newtheorem{example}[thm]{\protect\examplename}
  \theoremstyle{remark}
  \newtheorem{rem}[thm]{\protect\remarkname}

\usepackage{thmtools}

\usepackage{amsfonts}
\usepackage{german}
\def\ps@pprintTitle{%
 \let\@oddhead\@empty
 \let\@evenhead\@empty
 \def\@oddfoot{\centerline{\thepage}}%
 \let\@evenfoot\@oddfoot}

\makeatother

\usepackage{babel}
  \providecommand{\assumptionname}{Assumption}
  \providecommand{\corollaryname}{Corollary}
  \providecommand{\examplename}{Example}
  \providecommand{\remarkname}{Remark}
\providecommand{\theoremname}{Theorem}

\begin{document}
\global\long\def\E{\mathbb{E}}

\global\long\def\Er{\mathbb{\mathbf{Er}}}

\global\long\def\I{\mathbf{1}}

\global\long\def\N{\mathbb{N}}

\global\long\def\R{\mathbb{R}}

\global\long\def\C{\mathbb{C}}

\global\long\def\Q{\mathbb{Q}}

\global\long\def\P{\mathbb{P}}

\global\long\def\D{\mathcal{D}}

\global\long\def\dom{\operatorname{dom}}

\global\long\def\b#1{\mathbb{#1}}

\global\long\def\c#1{\mathcal{#1}}

\global\long\def\s#1{{\scriptstyle #1}}

\global\long\def\u#1#2{\underset{#2}{\underbrace{#1}}}

\global\long\def\r#1{\xrightarrow{#1}}

\global\long\def\mr#1{\mathrel{\raisebox{-2pt}{\ensuremath{\xrightarrow{#1}}}}}

\global\long\def\t#1{\left.#1\right|}

\global\long\def\f#1{\lfloor#1\rfloor}

\global\long\def\sc#1#2{\langle#1,#2\rangle}

\global\long\def\abs#1{\lvert#1\rvert}

\global\long\def\bnorm#1{\Bigl\lVert#1\Bigr\rVert}

\global\long\def\wraum{(\Omega,\c F,\P)}

\global\long\def\fwraum{(\Omega,\c F,\P,(\c F_{t}))}

\global\long\def\norm#1{\lVert#1\rVert}

\begin{frontmatter}{}

\title{Estimation error for occupation time functionals of stationary Markov
processes}

\author{Randolf Altmeyer}

\author{Jakub Chorowski}

\address{Institut für Mathematik, Humboldt-Universität zu Berlin, Unter den
Linden 6, 10099 Berlin, Germany}
\begin{abstract}
The approximation of integral functionals with respect to a stationary
Markov process by a Riemann-sum estimator is studied. Stationarity
and the functional calculus of the infinitesimal generator of the
process are used to get a better understanding of the estimation error
and to prove a general error bound. The presented approach admits
general integrands and gives a unifying explanation for different
rates obtained in the literature. Several examples demonstrate how
the general bound can be related to well-known function spaces. 
\end{abstract}

\end{frontmatter}{}

\section{Introduction}

Statistics for continuous-time Markov processes is usually based on
the observation of a sample path. Typically, only discrete-time observations
are available. An important task is the estimation of integral functionals
such as 
\[
\Gamma_{T}\left(f\right)=\int_{0}^{T}f\left(X_{r}\right)dr,\,\,\,\,T\geq0.
\]
Here, $X=(X_{r})_{r\geq0}$ is an $\c S$-valued Markov process on
some probability space $(\Omega,\c F,\P)$ for a Polish space $\c S$
equipped with its Borel-$\sigma$-field and $f:\c S\rightarrow\R$
is a given function such that $\Gamma_{T}(f)$ is well-defined. Functional
operators of this form appear in many problems. For instance, in mathematical
finance they are used to model path dependent derivatives (see \citet{hugonnier1999feynman},
\citet{Chesney1997}). In evolutionary dynamics the value $f(x)$
is often associated with the cost of staying in state $x$ (\citet{Pollett2003}).
The most important case for applications is $f=\I_{A}$ for a Borel
set $A$. $\Gamma_{T}(f)$ is then known as the \emph{occupation time}
of $X$ in $A$ and measures the time $X$ spends in $A$ before $T$.
For general $f$ and when $X$ is ergodic with invariant measure $\mu$,
integral functionals are also important for studying the long term
behavior of the process as $T\rightarrow\infty$, since $T^{-1}\Gamma_{T}(f)\rightarrow\int fd\mu$
by the ergodic theorem. Furthermore, the smoothness properties of
$x\mapsto\int_{0}^{T}f(x+X_{r})dr$ play an important role for solving
ordinary differential equations, for example in combination with the
phenomenon of regularization by noise (\citet{Catellier2016}). 

Our goal is to approximate $\Gamma_{T}(f)$ given the equidistant
observations $X_{k\triangle_{n}}$, $k=0,\dots,n-1$, where $\triangle_{n}=T/n$.
A natural candidate for this is the Riemann-sum estimator
\[
\hat{\Gamma}_{T,n}\left(f\right)=\sum_{k=1}^{n}f(X_{(k-1)\triangle_{n}})\triangle_{n}.
\]
Consistency of this estimator follows from Riemann-approximation already
under weak assumptions on $f$ and $X$. The rate of convergence,
however, depends on $f$, $X$, $T$ and $\triangle_{n}$. For deriving
useful finite sample bounds it is crucial to make these dependencies
explicit. The Riemann-sum estimator has appeared in many places in
the literature, mostly for estimating the occupation time (\citet{Chorowski2015a},
\citet{Gobet2016}), or as a proxy for approximating the local time
of $X$ in the diffusion case, such as in \citet{Hoffmann1999}. For
general $f$ see also \citet{Dion2016}. Error bounds are usually
derived ad-hoc, leading to suboptimal bounds or without explicit constants
for the dependence on parameters. 

It is clear that the approximation error depends on the smoothness
of $f$. Note, however, when measuring the approximation error in
the $L^{2}(\P)$ sense, it is not possible, in general, to have a
faster rate of convergence than $\sqrt{T}\triangle_{n}$, even with
smooth $f$. Indeed, for $X$ being a Brownian motion and $f$ the
identity, we have $\norm{\Gamma_{T}\left(f\right)-\hat{\Gamma}_{T,n}\left(f\right)}_{L^{2}\left(\P\right)}=\sqrt{T}\triangle_{n}/\sqrt{3}$.
A systematic study of this approximation problem has started only
recently. For one-dimensional diffusion processes with smooth coefficients
and $T=1$ \citet{Ngo2011} obtain the rate $n^{-3/4}$ for the special
case of indicator functions $f=\I_{[K,\infty)}$, $K\in\R$. Interestingly,
they also provide a lower bound for the rate $n^{-3/4}$ in the $L^{2}(\P)$-sense.
The specific analysis for indicator functions, however, cannot explain
which rates we can expect for more general $f$. For Markov processes
in $\R^{d}$ with transition kernels satisfying certain heat kernel
bounds \citet{Ganychenko2015} prove for bounded $f$ the $L^{p}(\P)$-bound
$\norm f_{\infty}\sqrt{T}(\triangle_{n}\log n)^{1/2}$, $p\geq2$.
$\alpha$-Hölder functions for $0<\alpha\leq1$ are studied by \citet{Kohatsu-Higa2014},
again for one-dimensional smooth diffusions, and by \citet{Ganychenko2015a},
for the same class of Markov processes as for \citet{Ganychenko2015}.
They essentially obtain the same upper bounds, $n^{-(1+\alpha)/2}$
for \citet{Kohatsu-Higa2014}, and $\norm f_{\alpha}\sqrt{T}\triangle_{n}^{(1+\alpha)/2}$
for \citet{Ganychenko2015a}, both of them losing a $\log n$-factor
when $\alpha=1$. Surprisingly, indicators obey the same bounds as
$1/2$-Hölder functions.

The aim of this paper is twofold. First, we want to study the approximation
of $\Gamma_{T}(f)$ by $\hat{\Gamma}_{T,n}(f)$ for more general functions
$f$ on arbitrary state spaces $\c S$ and to find a unifying mathematical
explanation for the different rates obtained in the literature. Second,
we want to develop a better understanding of the key quantities driving
the estimation error. For this, we focus on the special but important
case of stationary Markov processes, because this allows us to calculate
the error explicitly in terms of the associated semigroup. The main
insight of our results is that the discretization error depends on
the action of fractional powers of the infinitesimal generator applied
to $f$. We demonstrate in several examples how this can be related
to more familiar $L^{2}$-Sobolev norms of $f$. These norms are the
key idea to explain both the rates for Hölder and indicator functions
by suitable interpolation. The dependence on $T$ in the error is
explicit. This allows us, when $X$ is ergodic, to approximate integral
functionals with respect to the invariant measure of $X$ under weaker
conditions than the ones commonly used in the literature. Our approach
is based on the functional calculus of the generator. We therefore
consider only stationary Markov processes whose generators are normal
operators.

The paper is organized as follows. In Section \ref{sec:A-general-upper}
we state a general upper bound for the $L^{2}(\P)$-approximation
error. We apply it to approximate integral functionals with respect
to the invariant measure, when $X$ is ergodic. In Section \ref{sec:Examples}
we study several concrete examples of processes and functions. We
also discuss the important example of Brownian motion, which is not
a stationary Markov process, but which can be approximated by reflected
stationary diffusions. Proofs can be found in Section \ref{sec:Proofs}.
The appendix contains a brief summary of the most important facts
about semigroups and the functional calculus for normal operators.

\section{\label{sec:A-general-upper}A general upper bound}

In the following, we assume that $X$ is a continuous-time Markov
process on some probability space $(\Omega,\c F,\P)$ and with Polish
state space $\c S$. For any measure $\mu$ on $\c S$ we denote by
$L^{2}(\mu):=L^{2}(\c S,\mu)$ the space of square integrable functions
$f:\c S\rightarrow\R$ with respect to $\mu$ and with norm $\norm f_{\mu}=(\int f^{2}d\mu)^{1/2}$.
$\norm{\cdot}_{\infty,\mu}$ denotes the sup-norm in $L^{\infty}(\c{\mu})$
and $a\lesssim b$ for $a,b\in\R$ means $a\leq Cb$ for some constant
$C>0$. $Z_{n}=O_{\P}(a_{n})$ for a sequence of random variables
$(Z_{n})_{n\geq1}$ and real numbers $(a_{n})_{n\geq1}$ means that
$a_{n}^{-1}Z_{n}$ is tight. For basic concepts of semigroup theory
and functional calculus refer to the appendix. Our main assumptions
are the following:
\begin{assumption}
\label{assu:OT:Markov process assumption}$X$ is a stationary Markov
process with invariant measure $\mu$. The associated semigroup $(P_{r})_{r\geq0}$
is Feller and its infinitesimal generator $L$ is a normal operator
on $L^{2}(\mu).$ 
\end{assumption}
These assumptions are satisfied for many important processes. A leading
example is the standard Ornstein-Uhlenbeck process. A process with
normal but not necessarily self-adjoint generator will be discussed
in Section \ref{subsec:Infinite-dimensional-diffusions}. Observe
that for a stationary Markov process and $f\in L^{2}(\mu)$ both $\Gamma_{T}(f)$
and $\hat{\Gamma}_{T,n}(f)$ are $\mu$-a.s.-well-defined random variables
in $L^{2}(\P)$. We consider the spaces $\c D^{s}(L):=\text{dom}(|L|^{s/2})\subset L^{2}(\mu)$,
$s\geq0$, defined via the functional calculus of $L$, and with seminorm
$\norm f_{\c D^{s}(L)}:=\norm{|L|^{s/2}f}_{\mu}$. When $X$ is the
Ornstein-Uhlenbeck process, then the related spaces $\text{dom}((I-L)^{s/2})\subset\c D^{s}(L)$
are known as Bessel potential spaces and play an important role in
Malliavin calculus (\citet{watanabe1984lectures}). We are now ready
to state the general upper bound.
\begin{thm}
\label{thm:OT:Rates for Bessel potential spaces}Let $X$ be a Markov
process satisfying Assumption \ref{assu:OT:Markov process assumption}.
There exists a universal constant $C$ such that for all $0\leq s\leq1$
and $f\in\c D^{s}(L)$
\[
\left\Vert \Gamma_{T}(f)-\hat{\Gamma}_{T,n}(f)\right\Vert _{L^{2}(\P)}\leq C\|f\|_{\c D^{s}(L)}\sqrt{T}\triangle_{n}^{\frac{1+s}{2}}.
\]
\end{thm}
The proof of this theorem is remarkably short. For $s=0$ we have
$\c D^{0}(L)=L^{2}(\mu)$ and the rate is $\sqrt{T}\triangle_{n}^{1/2}$
which slightly improves the results of \citet{Ganychenko2015}, removing
an additional $\sqrt{\log n}$. Since $\c D^{s}(L)\subset\c D^{1}(L)$
for $s\geq1$, the rate is never better than $\sqrt{T}\triangle_{n}$.
For $0<s<1$ the bound interpolates between the two extreme cases.
A deeper understanding of the spaces $\c D^{s}(L)$ requires more
explicit knowledge about the generator. For example, if $L$ is self-adjoint,
then $|L|=-L$ and thus $\norm f_{\c D^{1}(L)}^{2}=\norm{(-L)^{1/2}f}_{\mu}^{2}=\left\langle -Lf,f\right\rangle _{\mu}$.
This is the Dirichlet form associated with $L$ and $\mu$. It is
typically easier to analyze than studying $\D^{1}(L)$ directly in
terms of the functional calculus. An important example are diffusions
on $\R^{d}$ such that for sufficiently smooth functions $f$ the
Dirichlet form is bounded by $\norm{\left|\nabla f\right|}_{\mu}^{2}$,
the $L^{2}(\mu)$-norm of the gradient of $f$. This immediately leads
to upper bounds for Hölder and indicator functions. Up to some additional
conditions, we will show that $\alpha$-Hölder functions lie in $\c D^{\alpha}(L)$
and indicator functions of certain cylinder sets of $\R^{d}$ lie
in $\c D^{1/2}(L)$. This gives a unifying explanation for the different
rates (see also Remark \ref{rem:}). These and other examples will
be discussed in Section \ref{sec:Examples}. A simple corollary shows
that the result of Theorem \ref{thm:OT:Rates for Bessel potential spaces}
remains valid if we relax the assumption of starting from the invariant
distribution.
\begin{cor}
\label{cor:Not stationary initial distribution}Let $X$ be a Markov
process satisfying Assumption \ref{assu:OT:Markov process assumption},
except that $X_{0}\overset{d}{\sim}\eta$ such that $\eta\ll\mu$
with density $d\eta/d\mu$. Then there exists a universal constant
$C$ such that for all $0\leq s\leq1$ and $f\in\c D^{s}(L)$
\[
\left\Vert \Gamma_{T}(f)-\hat{\Gamma}_{T,n}(f)\right\Vert _{L^{2}(\P)}\leq C\Big\|\frac{d\eta}{d\mu}\Big\|_{\infty,\mu}^{1/2}\|f\|_{\c D^{s}(L)}\sqrt{T}\triangle_{n}^{\frac{1+s}{2}}.
\]
\end{cor}
Instead of $\Gamma_{T}(f)$ a different target functional is often
$\int fd\mu$. It is well-known that $T^{-1}\Gamma_{T}(f)$ is $\sqrt{T}$-consistent
for $\int fd\mu$, i.e. $T^{-1}\Gamma_{T}(f)-\int fd\mu=O_{\P}(T^{-1/2})$,
when $L$ is self-adjoint and $f\in\text{dom}(L^{-1/2})$ (see e.g.
\citet{Kipnis1986}). By Theorem \ref{thm:OT:Rates for Bessel potential spaces}
we can now extend this to the estimator $T^{-1}\hat{\Gamma}_{T,n}(f)$
and more general $L^{2}(\mu)$-functions.
\begin{thm}
\label{thm:Convergence of the occupation measure to the stationary measure}Let
$X$ be a Markov process satisfying Assumption \ref{assu:OT:Markov process assumption}.
There exists a universal constant $C$ such that for all $f\in L^{2}(\mu)$
with $f_{0}\in\text{dom}(|L|^{-1/2})$, $f_{0}=f-\int fd\mu$,
\[
\left\Vert T^{-1}\hat{\Gamma}_{T,n}\left(f\right)-\int_{\mathcal{S}}f(x)\,d\mu(x)\right\Vert _{L^{2}(\P)}\leq\frac{C}{\sqrt{T}}\left(\|f\|_{\mu}\triangle_{n}^{1/2}+\left\Vert |L|^{-1/2}f_{0}\right\Vert _{\mu}\right).
\]
\end{thm}
As an example for $\text{dom}(|L|^{-1/2})$ being non-trivial, assume
that $0$ is a simple eigenvalue of $L$ and that $L$ has a spectral
gap, i.e. $s_{0}>0$, where $s_{0}=\sup\{r>0:B(0,r)\cap\sigma(L)=\{0\}\}$
and $B(0,r)=\{z\in\C:|z|\leq r\}$. In that case $X$ is ergodic and
it can be shown that $f_{0}\in\text{dom}(|L|^{-1/2})$ is satisfied
whenever $f$ is non-constant (\citet{bakry2013analysis}). Furthermore,
the upper bound of the theorem simplifies, since
\[
\left\Vert |L|^{-1/2}f_{0}\right\Vert _{\mu}\leq s_{0}^{-1/2}\|f_{0}\|_{\mu}\leq s_{0}^{-1/2}\|f\|_{\mu}.
\]
A concrete example of a process with spectral gap is the Ornstein-Uhlenbeck
process (Chapter 4 of \citet{bakry2013analysis}). In general, Theorem
\ref{thm:Convergence of the occupation measure to the stationary measure}
shows that in order to achieve the rate $\sqrt{T}$ as $n,T\rightarrow\infty$
there is essentially no gain in the high-frequency case, i.e. $\triangle_{n}\rightarrow0$,
compared to the low-frequency case with $\triangle_{n}$ fixed. The
error bound improves on the commonly used condition in the literature
that $T\triangle_{n}\lesssim1$ to achieve $\sqrt{T}$-consistency
(see e.g. Section 5 of \citet{Dion2016}). The theorem also shows
that the $T\rightarrow\infty$ case is controlled by negative powers
of the generator, while the fixed $T$ case is controlled by positive
powers of $|L|$ according to Theorem \ref{thm:OT:Rates for Bessel potential spaces}. 

\section{\label{sec:Examples}Examples}

In this section, we apply the general bound from Theorem \ref{thm:OT:Rates for Bessel potential spaces}
to several important examples. Our goal is not to give an exhaustive
list, but rather to demonstrate the theory in straightforward cases.
We first study Markov jump processes, i.e. continuous time Markov
processes with countable state spaces. We then consider a special
class of diffusion processes for which the spaces $\c D^{s}(L)$ can
be described easily via the Dirichlet form $\left\langle -Lf,f\right\rangle _{\mu}$.
After that, we show for the one-dimensional Brownian motion how the
assumption of stationarity can be removed. Finally, we discuss two
infinite dimensional diffusions, where the second one is an example
of a process whose infinitesimal generator is only a normal operator.

\subsection{Markov-jump processes}

Consider a continuous time Markov process $(X_{r})_{r\geq0}$ on a
countable state space $\mathcal{S}$. Such a process can always be
realized as $X_{r}=Y_{N_{r}}$ for a Markov chain $(Y_{s})_{s\in\c S}$
starting in some initial distribution $\mu$ with transition probabilities
$(P_{xy})_{x,y\in S}$ and an independent Poisson process $(N_{r})_{r\geq0}$
with intensity $0<\lambda<\infty$ (Chapter 4.2 of \citet{Ethier1986}).
Observing a path of the process $X$ at discrete times $0,\triangle_{n},2\triangle_{n},...,(n-1)\triangle_{n}$,
we can identify the jump times with $\triangle_{n}$ precision. Hence,
if the function $f$ is bounded, then every jump contributes at most
$2\|f\|_{\infty}$ to the estimation error $|\Gamma_{T}(f)-\hat{\Gamma}_{T,n}(f)|$
which yields the bound
\[
\norm{\Gamma_{T}(f)-\hat{\Gamma}_{T,n}(f)}_{L^{2}(\P)}\leq2\|f\|_{\infty}\E\left[N_{T}^{2}\right]^{1/2}\triangle_{n}=2\|f\|_{\infty}\left(\lambda T+\left(\lambda T\right)^{2}\right)^{1/2}\triangle_{n}.
\]
This already provides the optimal rate $\triangle_{n}$ but requires
the function $f$ to be bounded. Moreover, the error grows linearly
in $T$ as opposed to $\sqrt{T}$ as in Theorem \ref{thm:OT:Rates for Bessel potential spaces}.
We can improve on this assuming in addition that $X$ is stationary
with invariant measure $\mu$ and, for instance, reversible, i.e.
$P^{\top}=P$. The infinitesimal generator is $L=\lambda(P-I)$ (\citet{Ethier1986})
which is a bounded, non-negative self-adjoint operator. Therefore
$\norm f_{D^{1}(L)}\le\norm{(-L)^{1/2}}\norm f_{\mu}\leq\sqrt{\lambda}\norm f_{\mu}$
with operator norm $\norm{(-L)^{1/2}}$. It follows that $\c D^{1}(L)=L^{2}(\mu)$
and we obtain by Theorem \ref{thm:OT:Rates for Bessel potential spaces}
\[
\left\Vert \Gamma_{T}(f)-\hat{\Gamma}_{T,n}(f)\right\Vert _{L^{2}(\P)}\leq C\sqrt{\lambda}\|f\|_{\mu}\sqrt{T}\triangle_{n}.
\]

Note that the results of \citet{Ganychenko2015,Ganychenko2015a} do
not apply here, because the state space is countable and therefore
heat kernel bounds are not available.

\subsection{\label{subsec:Diffusions-with-generator}Diffusions with generator
in divergence form}

Let $(X_{r})_{r\geq0}$ be a stationary diffusion with values in some
closed subset $U\subset\c S:=\R^{d}$ with invariant measure $\mu$
that has support in $U$. In case $U\subsetneq\R^{d}$ we think of
$\mu$ as a measure on $\R^{d}$ and embed the domain of the infinitesimal
generator $\text{dom}(L)\subset L^{2}(U,\mu)$ canonically into $L^{2}(\R^{d},\mu)$
by letting $Lf:=L\tilde{f}$ whenever $f|_{U}=\tilde{f}$ for $f\in L^{2}(\R^{d},\mu)$,
$\tilde{f}\in L^{2}(U,\mu)$. We assume that $L$ is an elliptic operator
in divergence form (c.f. \citet[Chapter VII]{bass2006diffusions})
that satisfies 
\begin{equation}
\langle-Lf,g\rangle_{\mu}=\int_{\R^{d}}\left\langle A\left(x\right)\nabla f\left(x\right),\nabla g\left(x\right)\right\rangle _{\R^{d}}\,d\mu(x),\,\,\,\,f,g\in\text{dom}(L)\cap C^{2}(\R^{d}),\label{eq:OT:divergance form}
\end{equation}
for a measurable function $x\mapsto A(x)\in\R^{d\times d}$ such that
$A(x)$ is symmetric, positive definite for all $x\in\R^{d}$ and
such that $\norm{\left|A\right|}_{\infty,\mu}$ is finite, where $\left|\cdot\right|$
is any matrix norm. Observe that the right hand side is also well-defined
for $L^{2}(\R^{d},\mu)$-integrable functions $f,g\in C^{1}(\R^{d})$.
An operator $L$ satisfying \eqref{eq:OT:divergance form} is self-adjoint
on $\text{dom}(L)\cap C^{2}(\R^{d})$ such that $|L|^{1/2}=(-L){}^{1/2}$.
The starting point of our discussion is the observation that
\begin{align}
\norm f_{\c D^{s}(L)}^{2} & =\|(-L)^{s/2}f\|_{\mu}^{2}\leq\|(I-L)^{s/2}f\|_{\mu}^{2}\nonumber \\
 & \leq\|(I-L)^{1/2}f\|_{\mu}^{2}=\langle I-Lf,f\rangle_{\mu}\nonumber \\
 & \leq\norm f_{\mu}^{2}+\norm{\left|A\right|}_{\infty,\mu}\norm{\left|\nabla f\right|}_{\mu}^{2}\nonumber \\
 & \leq\max\left(1,\norm{\left|A\right|}_{\infty,\mu}\right)\norm f_{H^{1}(\mu)}^{2}\label{eq:upperBoundDivergenceForm-1}
\end{align}
for $f\in\text{dom}(L)\cap C^{2}(\R^{d})\subset\text{dom}(L)\subset\text{dom}((-L)^{s/2})=\c D^{s}(L)$
and $0\leq s\leq1$. The norm
\[
\norm f_{H^{1}(\mu)}=\norm f_{\mu}+\norm{\left|\nabla f\right|}_{\mu}
\]
is the $\mu$-weighted Sobolev norm. Combining this with Theorem \ref{thm:OT:Rates for Bessel potential spaces}
yields 
\begin{align}
 & \left\Vert \Gamma_{T}(f)-\hat{\Gamma}_{T,n}(f)\right\Vert _{L^{2}(\P)}\nonumber \\
 & \leq\begin{cases}
C\max\left(1,\norm{\left|A\right|}_{\infty,\mu}^{1/2}\right)\|f\|_{H^{1}(\mu)}\sqrt{T}\triangle_{n}, & f\in\text{dom}(L)\cap C^{2}(\R^{d}),\\
C\|f\|_{\mu}\sqrt{T}\triangle_{n}^{1/2}, & f\in L^{2}(\mu).
\end{cases}\label{eq:boundForH1_mu}
\end{align}
By interpolation between the two cases $f\in L^{2}(\mu)$ and $f\in\text{dom}(L)\cap C^{2}(\R^{d})$
we will study Hölder and indicator functions, recovering previously
obtained results. We will further give a unifying view on these results
by using Sobolev spaces instead. Before doing this let us discuss
some important examples where \eqref{eq:OT:divergance form} holds. 
\begin{example}[Ornstein-Uhlenbeck process]
\label{ex:OT:Ornstein-Uhlenbeck-process}Assume that $(X_{r})_{r\geq0}$
satisfies the stochastic differential equation 
\[
dX_{r}=-X_{r}dr+\sqrt{2}dW_{r}
\]
in $\R^{d}$ where $(W_{r})_{r\geq0}$ is a $d$-dimensional Brownian
motion. If $X_{0}\overset{d}{\sim}\mu$, where $\mu$ has Lebesgue
density $d\mu(x)/dx=(2\pi)^{-d/2}\exp(-|x|^{2}/2)$, then $X$ is
stationary with invariant measure $\mu$. The infinitesimal generator
$L$ satisfies
\begin{equation}
Lf(x)=-\left\langle x,\nabla f\left(x\right)\right\rangle _{\R^{d}}+\Delta f\left(x\right),\,\,\,\,x\in\R^{d},\label{eq:genOfOU}
\end{equation}
with $f\in\text{dom}(L)=H^{2}(\mu)$, the $\mu$-weighted Sobolev
space of twice weak differentiable functions with all partial derivatives
up to order two belonging to $L^{2}(\mu)$ (\citet{Chojnowska-Michalik2002}).
Using integration by parts it follows (c.f. \citet{pavliotis2014stochastic},
Section 4.4) that
\begin{equation}
\langle-Lf,g\rangle_{\mu}=\int_{\R^{d}}\left\langle \nabla f\left(x\right),\nabla g\left(x\right)\right\rangle _{\R^{d}}\,d\mu(x),\,\,\,\,f,g\in C^{2}(\R^{d}),.\label{eq:OT: D^1 for OU}
\end{equation}
Hence $L$ is a self-adjoint operator of the form \eqref{eq:OT:divergance form}
with $A_{jk}=\I(j=k)$ for all $1\leq j,k\leq d$. This example can
be generalized considerably (see \citet{Chojnowska-Michalik2002}
and Subsection \ref{subsec:Infinite-dimensional-diffusions} below). 
\end{example}
\begin{example}[Scalar diffusion with possibly attracting boundaries]
\label{ex:OT:Stationary-scalar-diffusion}Fix boundaries $-\infty\leq\lambda<\rho\leq\infty.$
Assume that $(X_{r})_{r\geq0}$ is a stationary diffusion process
on $[\lambda,\rho]$ solving the one-dimensional stochastic differential
equation 
\begin{equation}
dX_{r}=b(X_{r})dr+\sigma(X_{r})dW_{r},\label{eq:OT:SDE for scalar diffusion}
\end{equation}
for a continuous drift $b:[\lambda,\rho]\to\R$, strictly positive
continuous volatility $\sigma:[\lambda,\rho]\to(0,\infty)$ and a
one-dimensional Brownian motion $(W_{r})_{r\geq0}$. Sufficient conditions
for the existence of such a process can be found in \citet{Hansen1998}.
In particular, stationarity is guaranteed if the \emph{speed density}
\[
m\left(x\right)=\frac{1}{\sigma^{2}\left(x\right)}\exp\Big(\int_{x_{0}}^{x}\frac{2b(y)}{\sigma^{2}(y)}dy\Big),\,\,\,\,\lambda\leq x_{0}\leq\rho,\,\lambda<x<\rho,
\]
is integrable on $[\lambda,\rho]$. The stationary measure then has
the density
\[
\frac{d\mu\left(x\right)}{dx}=C_{0}m\left(x\right)\I\left(\lambda<x<\rho\right),
\]
where $C_{0}$ is a normalizing constant. The infinitesimal generator
$L$ satisfies
\begin{align*}
Lf\left(x\right) & =b(x)f'(x)+\frac{\sigma^{2}(x)}{2}f''(x)\\
 & =\frac{1}{2}\left(\frac{d\mu(x)}{dx}\right)^{-1}\left(f'(x)\sigma^{2}(x)\frac{d\mu(x)}{dx}\right)',\,\,\,\,\lambda<x<\rho,
\end{align*}
with $f\in\text{dom}(L)$, where 
\begin{align*}
\operatorname{dom}(L) & =\Big\{ f\in L^{2}([\lambda,\rho],\mu):f\text{ and\,}\,f'\text{ are absolutely continuous with}\\
 & \qquad\lim_{x\searrow\lambda}f'(x)m\left(x\right)\sigma^{2}\left(x\right)=\lim_{x\nearrow\rho}f'(x)m\left(x\right)\sigma^{2}\left(x\right)=0\text{ and }\\
 & \qquad Lf\in L^{2}([\lambda,\rho],\mu)\Big\}.
\end{align*}
For details see Section 3.3 of \citet{Hansen1998}. Embedding the
domain into $L^{2}(\mu)$ as mentioned before and integrating by parts
it follows that
\begin{equation}
\langle-Lf,g\rangle_{\mu}=\int_{\R}f'(x)g'(x)\sigma^{2}(x)\,d\mu(x),\,\,\,\,f,g\in\text{dom}(L)\cap C^{2}(\R),\label{eq:divergenceFormForScalarDiffusion}
\end{equation}
which is of the form \eqref{eq:OT:divergance form} with $A=\sigma^{2}$.
For $b(x)=-x$ and $\sigma(x)=\sqrt{2}$, $X$ is just the one-dimensional
Ornstein-Uhlenbeck process. 
\end{example}
\begin{example}[Reflected diffusion]
\label{ex:OT:Scalar-diffusion-process with reflection}Fix boundaries
$-\infty<\lambda<\rho<\infty.$ Assume that $X$ is a one-dimensional
reflected diffusion on $[\lambda,\rho]$. By this we mean that $X$
satisfies the Skorokhod type stochastic differential equation
\begin{equation}
dX_{r}=b(X_{r})dr+\sigma(X_{r})dW_{r}+dK_{r},\label{eq:reflectedDiffusion}
\end{equation}
for a bounded measurable drift $b:[\lambda,\rho]\to\R$, strictly
positive continuous volatility $\sigma:[\lambda,\rho]\to(0,\infty)$,
$(W_{r})_{r\geq0}$ is a Brownian motion and $(K_{r})_{r\geq0}$ is
an adapted continuous process with finite variation starting from
$0$ and such that for every $r\geq0$ 
\[
\int_{0}^{r}\I_{(\lambda,\rho)}(X_{s})dK_{s}=0.
\]
The invariant measure and the generator $L$ are as in the last example.
Since $[\lambda,\rho]$ is compact, the domain simplifies to
\begin{align*}
\operatorname{dom}(L) & =\Big\{ f\in L^{2}([\lambda,\rho],\mu):f\text{ and\,}\,f'\text{ are absolutely continuous with}\\
 & \qquad f'(\lambda)=f'(\rho)=0\text{ and }Lf\in L^{2}([\lambda,\rho],\mu)\Big\}.
\end{align*}
Therefore \eqref{eq:OT:divergance form} holds here, as well. For
more details see \citet{Chorowski2015a}.
\end{example}

\subsubsection{Hölder functions}

Consider an $\alpha$-Hölder continuous function $f:\R^{d}\rightarrow\R$,
i.e. $f$ has finite Hölder-norm
\[
\norm f_{\alpha}=\sup_{x\neq y\in\R^{d}}\frac{|f(x)-f(y)|}{|x-y|^{\alpha}},
\]
$0\leq\alpha\leq1$, and such that $f\in L^{2}(\mu)$. We want to
derive an upper bound for $\norm{\Gamma_{T}(f)-\hat{\Gamma}_{T,n}(f)}_{L^{2}(\P)}$
in terms of $\norm f_{\alpha}$ using \eqref{eq:boundForH1_mu}. Let
$(\varphi_{\varepsilon})_{\varepsilon\geq0}$ be a non-negative smooth
kernel, i.e. $\varphi_{\varepsilon}(x)=\varepsilon^{-1}\varphi(\varepsilon^{-1}x)$,
$0\leq\varphi\in C_{c}^{\infty}(\R^{d})$, $\text{supp}(\varphi)\subset[-1,1]^{d}$,
$\int_{\R^{d}}\varphi(x)\,dx=1$. Then the convolution $f_{\varepsilon}=f*\varphi_{\varepsilon}$
lies in $C^{\infty}(\R^{d})$ with bounded derivatives and therefore
$f_{\varepsilon}\in L^{2}(\mu)\cap C^{2}(\R^{d})$. In order to apply
\eqref{eq:boundForH1_mu} to $f_{\varepsilon}$ we also need that
$f_{\varepsilon}\in\text{dom}(L)$. From the examples above we see
that the only issue might be possible boundary conditions for functions
in $\text{dom}(L)$. In order to remove these conditions, we assume
the following:
\begin{assumption}
\label{assu:densityAssumptionForGenerator}$\text{dom}(L)\cap C^{2}(\R^{d})$
is dense in $L^{2}(\mu)\cap C^{1}(\R^{d})$ with respect to $\norm{\cdot}_{H^{1}(\mu)}$. 
\end{assumption}
This assumption is relatively weak. It is satisfied in all examples
above. In particular, if functions $f\in\text{dom}(L)$ do not have
to satisfy any boundary conditions, then $L^{2}(\mu)\cap C^{2}(\R^{d})\subset\text{dom}(L)\cap C^{2}(\R^{d})$,
as is the case for the Ornstein-Uhlenbeck process. By approximation
we can thus extend \eqref{eq:boundForH1_mu} to 
\begin{align}
 & \left\Vert \Gamma_{T}(f)-\hat{\Gamma}_{T,n}(f)\right\Vert _{L^{2}(\P)}\nonumber \\
 & \leq\begin{cases}
C\max\left(1,\norm{\left|A\right|}_{\infty,\mu}^{1/2}\right)\|f\|_{H^{1}(\mu)}\sqrt{T}\triangle_{n}, & f\in L^{2}(\mu)\cap C^{1}(\R^{d}),\\
C\|f\|_{\mu}\sqrt{T}\triangle_{n}^{1/2}, & f\in L^{2}(\mu).
\end{cases}\label{eq:boundForH1_mu-1}
\end{align}

We will extend this to Sobolev functions in Subsection \ref{subsec:Sobolev-functions}.
It follows that $f_{\varepsilon}\in L^{2}(\mu)\cap C^{1}(\R^{d})$.
Using $\int\varphi\left(x\right)dx=1$ and $\int\nabla\varphi(x)dx=0$,
we have 
\begin{align*}
\norm{f-f_{\varepsilon}}_{\mu}^{2} & =\int\left|\int\left(f\left(x\right)-f\left(x+\varepsilon y\right)\right)\varphi\left(y\right)dy\right|^{2}d\mu\left(x\right)\leq\norm f_{\alpha}^{2}\varepsilon^{2\alpha},\\
\norm{\left|\nabla f_{\varepsilon}\right|}_{\mu}^{2} & =\int\left|\frac{f\left(x\right)}{\varepsilon}\int\nabla\varphi\left(y\right)dy-\nabla f_{\varepsilon}\left(x\right)\right|^{2}d\mu\left(x\right)\\
 & =\frac{1}{\varepsilon^{2}}\int\left|\int\left(f\left(x\right)-f\left(x+\varepsilon y\right)\right)\nabla\varphi\left(y\right)dy\right|^{2}d\mu\left(x\right)\\
 & \lesssim\norm f_{\alpha}^{2}\varepsilon^{2\alpha-2}.
\end{align*}
Hence $\norm{f_{\varepsilon}}_{H^{1}(\mu)}\lesssim\norm f_{\mu}+\norm f_{\alpha}\varepsilon^{\alpha-1}$.
Together with $\norm{f_{\varepsilon}}_{\mu}\leq\norm{f-f_{\varepsilon}}_{\mu}+\norm f_{\mu}$
and \eqref{eq:boundForH1_mu-1} this yields
\begin{align*}
 & \left\Vert \Gamma_{T}(f)-\hat{\Gamma}_{T,n}(f)\right\Vert _{L^{2}(\P)}\\
 & \qquad\leq\left\Vert \Gamma_{T}(f-f_{\varepsilon})-\hat{\Gamma}_{T,n}(f-f_{\varepsilon})\right\Vert _{L^{2}(\P)}+\left\Vert \Gamma_{T}(f_{\varepsilon})-\hat{\Gamma}_{T,n}(f_{\varepsilon})\right\Vert _{L^{2}(\P)}\\
 & \qquad\lesssim\norm f_{\alpha}\sqrt{T}\triangle_{n}^{1/2}\varepsilon^{\alpha}+\norm f_{\alpha}\sqrt{T}\triangle_{n}\varepsilon^{\alpha-1}+\norm f_{\mu}\sqrt{T}\triangle_{n}.
\end{align*}
Choosing $\varepsilon=\triangle_{n}^{1/2}$ implies the bound $\norm f_{\alpha}\sqrt{T}\triangle_{n}^{(1+\alpha)/2}+\norm f_{\mu}\sqrt{T}\triangle_{n}$.
Up to the second term, which is of smaller order as long as $\alpha<1$,
these are the rates obtained by \citet{Kohatsu-Higa2014} and \citet{Ganychenko2015a}.
If $L$ satisfies a Poincaré type inequality, i.e. if there exists
a constant $c<\infty$ such that for all $f\in\text{dom}(L)$
\begin{equation}
\norm{f_{0}}_{\mu}\leq c\norm{\left|\nabla f\right|}_{\mu},\label{eq:poincare}
\end{equation}
where $f_{0}=f-\int fd\mu$, then $\norm f_{\alpha}=\norm{f_{0}}_{\alpha}$
and thus
\begin{align}
 & \left\Vert \Gamma_{T}(f)-\hat{\Gamma}_{T,n}(f)\right\Vert _{L^{2}(\P)}=\left\Vert \Gamma_{T}(f_{0})-\hat{\Gamma}_{T,n}(f_{0})\right\Vert _{L^{2}(\P)}\nonumber \\
 & \qquad\lesssim\norm f_{\alpha}\sqrt{T}\triangle_{n}^{(1+\alpha)/2}.\label{eq:}
\end{align}
Poincaré inequalities hold for many stationary measures $\mu$, for
example for the Ornstein-Uhlenbeck process and in Example \ref{ex:OT:Stationary-scalar-diffusion}
when $m(x)$ is uniformly bounded from above and below. For other
examples see \citet{bakry2013analysis} and \citet{chen2006eigenvalues}.
Observe that for $\alpha=1$ the upper bound is $\norm f_{1}\sqrt{T}\triangle_{n}$,
removing an additional $\sqrt{\log n}$ term present in the results
of \citet{Kohatsu-Higa2014} and \citet{Ganychenko2015a}. We sum
up this discussion in the following theorem.
\begin{thm}
\label{thm:OT:Rate for Holder}Let $X$ be a stationary diffusion
with values in $\R^{d}$, invariant measure $\mu$ and whose generator
$L$ satisfies \eqref{eq:OT:divergance form} and Assumption \ref{assu:densityAssumptionForGenerator}.
Then we have for an $\alpha$-Hölder continuous function $f$, $0\leq\alpha\leq1$, 

\[
\left\Vert \Gamma_{T}(f)-\hat{\Gamma}_{T,n}(f)\right\Vert _{L^{2}(\P)}\lesssim\norm f_{\alpha}\sqrt{T}\triangle_{n}^{\frac{1+\alpha}{2}}+\norm f_{\mu}\sqrt{T}\triangle_{n}.
\]
If $L$ satisfies a Poincaré type inequality as in \eqref{eq:poincare}
for some $c<\infty$, then the upper bound is just $\norm f_{\alpha}\sqrt{T}\triangle_{n}^{\frac{1+\alpha}{2}}$.
\end{thm}

\subsubsection{Indicator functions}

We study now indicator functions $f=\I_{[K,\infty)}$, $K\in\R$,
such that $f\in L^{2}(\mu)$. We will argue again with \eqref{eq:boundForH1_mu-1}
by approximation. Let $(\varphi_{\varepsilon})_{\varepsilon>0}$ be
a non-negative smooth kernel as in the previous example, but this
time in $\R$. Then $f_{\varepsilon}=f*\varphi_{\varepsilon}$ is
bounded by $1$ and lies in $L^{2}(\mu)\cap C^{2}(\R)$. Moreover,
$f-f_{\varepsilon}$ has support in $[K-\varepsilon,K+\varepsilon]$
such that 
\begin{align*}
\norm{f-f_{\varepsilon}}_{\mu}^{2} & \leq\int_{K-\varepsilon}^{K+\varepsilon}d\mu,\\
\norm{f_{\varepsilon}'}_{\mu}^{2} & =\int\left|\frac{f\left(x\right)}{\varepsilon}\int\varphi'\left(y\right)dy-f_{\varepsilon}'\left(x\right)\right|^{2}d\mu\left(x\right)\\
 & =\frac{1}{\varepsilon^{2}}\int_{K-\varepsilon}^{K+\varepsilon}\left|\int\left(f\left(x\right)-f\left(x+\varepsilon y\right)\right)\varphi'\left(y\right)dy\right|^{2}d\mu\left(x\right)\\
 & \lesssim\frac{1}{\varepsilon^{2}}\int_{K-\varepsilon}^{K+\varepsilon}d\mu.
\end{align*}
As before $\norm{f_{\varepsilon}}_{\mu}\leq\norm{f-f_{\varepsilon}}_{\mu}+\norm f_{\mu}$.
If $\mu$ is absolutely continuous with respect to the Lebesgue measure
with bounded density, then $\varepsilon^{-1}\int_{K-\varepsilon}^{K+\varepsilon}d\mu$
is bounded and in that case we have from \eqref{eq:boundForH1_mu}
uniform in $K$
\begin{equation}
\left\Vert \Gamma_{T}(f)-\hat{\Gamma}_{T,n}(f)\right\Vert _{L^{2}(\P)}\lesssim\sqrt{T}\left(\triangle_{n}\varepsilon\right)^{1/2}+\sqrt{T}\triangle_{n}\varepsilon^{-1/2}+\sqrt{T}\triangle_{n}.\label{eq:indicatorBound}
\end{equation}
The last term is of lower order compared to the first two. Hence choosing
$\varepsilon=\triangle_{n}^{1/2}$ yields the rate $\sqrt{T}\triangle_{n}^{3/4}$
obtained by \citet{Ngo2011} (see also \citet{Kohatsu-Higa2014})
for one-dimensional diffusions. However, now the rate is uniform in
$K$ with explicit dependence on $T$. These arguments can easily
be extended to general dimensions and we obtain the following theorem.
\begin{thm}
\label{thm:OT:Rate for indicators}Let $X$ be a stationary diffusion
with values in $\R^{d}$, invariant measure $\mu$ and whose generator
$L$ satisfies \eqref{eq:OT:divergance form} and Assumption \ref{assu:densityAssumptionForGenerator}.
We assume furthermore that $\mu$ has bounded Lebesgue density. If
$f$ is an indicator function in $\R^{d}$ of the form $[K_{1},L_{1})\times\dots\times[K_{d},L_{d})$,
$-\infty<K_{j}<L_{j}\leq\infty$, $1\leq j\leq d$, then 
\[
\left\Vert \Gamma_{T}(f)-\hat{\Gamma}_{T,n}(f)\right\Vert _{L^{2}(\P)}\lesssim\sqrt{T}\triangle_{n}^{3/4},
\]
uniformly in $K_{j},L_{j}$.
\end{thm}
The same rate clearly holds up to constants for finite linear combinations
of such indicators.

\subsubsection{\label{subsec:Sobolev-functions}Sobolev functions}

Our goal is to explain why the bound for indicator functions in Theorem
\ref{thm:OT:Rate for indicators} can be obtained by evaluating the
bound for Hölder functions in Theorem \ref{thm:OT:Rate for Holder}
at $\alpha=1/2$. The key idea is to study indicators as elements
of certain Sobolev spaces and interpolate between the two extreme
cases in \eqref{eq:boundForH1_mu-1}. The closure of $L^{2}(\mu)\cap C^{1}(\R^{d})$
with respect to $\norm{\cdot}_{H^{1}(\mu)}$ yields the space $H^{1}(\mu)$,
a $\mu$-weighted Sobolev space. This is not a Banach space in general
(\citet{10.2307/40302922}). In order to avoid this issue we assume
that $\mu$ has a bounded Lebesgue density $d\mu/dx$. Then $L^{2}(\R^{d}):=L^{2}(\R^{d},\lambda)\subset L^{2}(\R^{d},\mu)$,
where $\lambda$ is the Lebesgue measure, and we have 
\[
\norm f_{H^{1}\left(\mu\right)}\leq\norm{\frac{d\mu}{dx}}_{\infty}\norm f_{H^{1}},\,\,\,\,f\in L^{2}(\R^{d})\cap C^{1}(\R^{d}).
\]
Here, $\norm{\cdot}_{\infty}:=\norm{\cdot}_{\infty,\lambda}$ is the
$\lambda$-sup-norm and $\norm f_{H^{1}}:=\norm f_{H^{1}(\lambda)}$
is the classical Sobolev norm with respect to $\lambda$. Taking the
closure of $L^{2}(\mu)\cap C^{2}(\R^{d})$ with respect to $\norm{\cdot}_{H^{1}}$
leads to the Sobolev space $H^{1}(\R^{d})$ of weakly differentiable
functions. This yields then instead of \eqref{eq:boundForH1_mu-1}
\begin{align}
 & \left\Vert \Gamma_{T}(f)-\hat{\Gamma}_{T,n}(f)\right\Vert _{L^{2}(\P)}\nonumber \\
 & \leq\begin{cases}
C\max\left(1,\norm{\left|A\right|}_{\infty,\mu}^{1/2}\right)\norm{\frac{d\mu}{dx}}_{\infty}^{1/2}\|f\|_{H^{1}}\sqrt{T}\triangle_{n}, & f\in H^{1}\left(\R^{d}\right),\\
C\norm{\frac{d\mu}{dx}}_{\infty}^{1/2}\|f\|_{\mu}\sqrt{T}\triangle_{n}^{1/2}, & f\in L^{2}(\R^{d}).
\end{cases}\label{eq:boundForH1}
\end{align}
Interpolation between the norms $\norm{\cdot}_{H^{1}}$ and $\norm{\cdot}_{\mu}$
is classical (\citet{adams2003sobolev}) and leads to the fractional
Sobolev spaces $H^{s}(\R^{d})$, $0\leq s\leq1$, with finite norm
\begin{align*}
\norm f_{H^{s}} & =\Big(\frac{1}{(2\pi)^{d/2}}\int_{\R^{d}}(1+|u|)^{2s}|\c Ff(u)|^{2}du\Big)^{1/2},
\end{align*}
where $\c Ff$ is the $L^{2}$-Fourier transform of $f$. Interpolation
theory for the operator $\Gamma_{T}(f)-\hat{\Gamma}_{T,n}(f)$ (Section
7.2 of \citet{adams2003sobolev}) yields for $f\in H^{s}(\R^{d})$
the upper bound 
\begin{equation}
\left\Vert \Gamma_{T}(f)-\hat{\Gamma}_{T,n}(f)\right\Vert _{L^{2}(\P)}\leq C\max\left(1,\norm{\left|A\right|}_{\infty,\mu}^{1/2}\right)\norm{\frac{d\mu}{dx}}_{\infty}^{1/2}\|f\|_{H^{s}}\sqrt{T}\triangle_{n}^{\frac{1+s}{2}}.\label{eq:SobolevRate}
\end{equation}

We sum up this discussion in a theorem.
\begin{thm}
\label{thm:OT:Rate for Sobolev functions}Let $X$ be a stationary
diffusion with values in $\R^{d}$, invariant measure $\mu$ and whose
generator $L$ satisfies \eqref{eq:OT:divergance form} and Assumption
\ref{assu:densityAssumptionForGenerator}. We assume furthermore that
$\mu$ has bounded Lebesgue density. Then \eqref{eq:SobolevRate}
holds for any $f\in H^{s}(\R^{d})$, $0\leq s\leq1$.
\end{thm}
Let us check what this bound implies for Hölder and indicator functions.
If $f$ is an $\alpha$-Hölder continuous function with $0<\alpha\leq1$
and compact support in $\R^{d}$, then it is well-known that $f\in H^{\alpha-\delta}(\R^{d})$
for any small $\delta>0$. \eqref{eq:SobolevRate} thus implies the
rate $\sqrt{T}\triangle_{n}^{(1+\alpha)/2-\delta}$. Using large deviation
bounds for $X$, we can extend this to general $\alpha$-Hölder continuous
$f$ losing a $\log\triangle_{n}$-factor in the rate. Consider now
the one-dimensional indicator $f=\I_{[K,L)}$, $-\infty<K<L<\infty$.
Then \eqref{eq:indicatorBound} remains true, but we also have $f\in H^{s}(\R)$
for all $s<1/2$ and \eqref{eq:SobolevRate} yields up to constants
the bound $\norm f_{H^{1/2-\delta}}\sqrt{T}\triangle_{n}^{3/4-\delta}$
for any small $\delta>0$. This is not uniform in $K,L$ anymore,
but still describes well how the error depends on $f$. Compared to
Theorems \ref{thm:OT:Rate for Holder} and \ref{thm:OT:Rate for indicators}
we see that explicit approximation via \eqref{eq:boundForH1_mu-1}
yields in general sharper bounds than \eqref{eq:SobolevRate}, but
the latter one is easier to apply since we only need to bound the
Fourier transform of $f$. 
\begin{rem}
\label{rem:}\begin{enumerate}[(i)] \item Assume that $H^{1}(\mu)$
is a Banach space. This is true, for instance, in the Examples \ref{ex:OT:Ornstein-Uhlenbeck-process}
and \ref{ex:OT:Stationary-scalar-diffusion}, when $m(x)$ is uniformly
bounded from above and below. In that case we can directly interpolate
between $H^{1}(\mu)$ and $L^{2}(\mu)$ with a similar bound as in
\eqref{eq:SobolevRate}, but with $\norm{\cdot}_{H^{s}}$ replaced
by an appropriate interpolation norm. The results in Theorems \ref{thm:OT:Rate for Holder}
and \ref{thm:OT:Rate for indicators} are explicit cases of this.
Up to boundary conditions this implies that $\alpha$-Hölder functions
lie in $\c D^{\alpha}(L)$, $0\leq\alpha\leq1$, and indicator functions
$f=\I_{[K,L)}$ lie in $\c D^{1/2}(L)$, $-\infty<K<L\leq\infty$.
\item Depending on the boundary conditions for functions $f\in\text{dom}(L)$
and if $\mu$ has bounded Lebesgue-density, in many examples it can
be shown that $H^{1}(\R^{d})$ embeds continuously into $\text{dom}((I-L)^{1/2})\subset\c D^{1}(L)$.
This holds, for instance, for the Ornstein-Uhlenbeck process and for
the reflected diffusions in Example \ref{eq:OT:SDE for scalar diffusion}.
Since $L^{2}(\R^{d})\subset\c D^{0}(L)=L^{2}(\R^{d},\mu)$, we obtain
by interpolation that $H^{s}(\R^{d})$ embeds continuously into $\text{dom}((I-L)^{s/2})\subset\c D^{s}(L)$.
In particular, the indicator functions $f=\I_{[K,L)}$ are elements
of $\c D^{1/2-\delta}(L)$ for any small $\delta>0$.\item Arguing
like in the proof of Corollary \ref{cor:Not stationary initial distribution}
the strict stationarity assumption can be relaxed and we obtain the
same results when $X$ starts from some initial distribution $\eta$
which is absolutely continuous with respect to $\mu$ and has bounded
density $d\eta/d\mu$. This is also true for Theorems \ref{thm:OT:Rate for Holder}
and \ref{thm:OT:Rate for indicators}. \end{enumerate}
\end{rem}

\subsection{Brownian motion and related diffusions}

We will demonstrate now how in some cases the procedure of the last
subsection can be adapted if the generator satisfies \eqref{eq:OT:divergance form},
but the invariant measure is not a probability measure. The main example
for this is the one-dimensional Brownian motion. Its generator is
$Lf=(1/2)\Delta f$, $\text{dom}(L)=H^{2}(\R)$, and the invariant
measure $\mu$ is the Lebesgue measure on $\R^{d}$. The key idea
idea is to approximate the process by stationary diffusions with reflecting
boundaries. For simplicity, we discuss only the one-dimensional case
under the additional assumption that the process starts from a distribution
$\eta$ which is absolutely continuous with respect to the Lebesgue
measure. 
\begin{thm}
\label{thm:BrownianMotion}Assume that $(X_{r})_{r\geq0}$ satisfies
$X_{r}=X_{0}+B_{r}$, where $(B_{r})_{r\geq0}$ is a one-dimensional
Brownian motion, $X_{0}\overset{d}{\sim}\eta$ is independent of $(B_{r})_{r\geq0}$
and $\eta$ has Lebesgue density with compact support. Then Theorem
\ref{thm:OT:Rate for Sobolev functions} applies to $(X_{r})_{r\geq0}$,
as well, with $A=1$ and $d\mu/dx=d\eta/dx$.
\end{thm}
Following a similar program as in \citet{Ngo2011} or \citet{Kohatsu-Higa2014}
this result can be extended to more general one-dimensional diffusion
processes, using the Lamperti and Girsanov transforms.

\subsection{\label{subsec:Infinite-dimensional-diffusions}Infinite dimensional
diffusions}

Since the general state space $\c S$ of $X$ is Polish, we can also
study infinite dimensional diffusions. Note that the results of \citet{Ganychenko2015,Ganychenko2015a}
do not apply here, because, in general, heat kernel bounds are not
available in this setting. Example \ref{ex:OT:Ornstein-Uhlenbeck-process}
can be generalized considerably. If $X$ satisfies the stochastic
differential equation
\[
dX_{r}=AX_{r}dr+Q^{1/2}dW_{r},
\]
where $A$ and $Q$ are operators on a separable Hilbert space $\c H$,
with $Q$ being bounded self-adjoint, then $X$ is a Gaussian Markov
process and the generator $L$ satisfies a similar formula as in \eqref{eq:genOfOU}
with $\nabla$ and $\Delta$ replaced by the corresponding Fréchet
derivatives $D$ and $D^{2}$. Under certain conditions on $A$ and
$Q$ the generator is reversible and $X$ has some invariant measure
$\mu$. The domain is again a $\mu$-weighted Sobolev space and the
associated Dirichlet form is
\[
\left\langle -Lf,g\right\rangle _{\mu}=\frac{1}{2}\int_{\c H}\left\langle Q^{1/2}Df\left(x\right),Q^{1/2}Dg\left(x\right)\right\rangle _{\c H}\,d\mu\left(x\right).
\]
The results of Section \ref{subsec:Diffusions-with-generator} therefore
remain formally the same. For details see \citet{Chojnowska-Michalik2002}.
A different kind of example are infinite dimensional systems of the
form 
\[
dX_{r}^{(i)}=\left(pV'\left(X_{r}^{(i+1)}-X_{r}^{(i)}\right)-qV'\left(X_{r}^{(i)}-X_{r}^{(i-1)}\right)\right)dr+dW_{r}^{(i)},
\]
$(r,i)\in[0,\infty)\times\mathbb{Z}$, where $\{(W_{r}^{(i)})_{r\geq0}:i\in\mathbb{Z}\}$
is an independent family of Brownian motions, $p,q\geq0$ with $p+q=1$
and $V$ is some potential function (\citet{Diehl2016a}). The authors
show stationarity of $X=(X_{r}^{(i)})_{r\geq0,i\in\mathbb{Z}}$, derive
the infinitesimal generator and derive its Dirichlet form (for the
symmetric part). A similar analysis can then be derived as in Section
\ref{subsec:Diffusions-with-generator}. Note that the generator in
this example is not necessarily self-adjoint, but it is always normal
(Lemma 2.1 of \citet{Diehl2016a}).

\section{\label{sec:Proofs}Proofs}

\subsection{Proof of Theorem \ref{thm:OT:Rates for Bessel potential spaces}}
\begin{proof}
We first assume that $f\in L^{2}(\mu)$. Expanding the squared error
yields
\begin{align*}
 & \norm{\Gamma_{T}\left(f\right)-\hat{\Gamma}_{T,n}\left(f\right)}_{L^{2}\left(\P\right)}^{2}=\E\Big[\Big|\sum_{k=1}^{n}\int_{\left(k-1\right)\triangle_{n}}^{k\triangle_{n}}\left(f\left(X_{r}\right)-f\left(X_{(k-1)\triangle_{n}}\right)\right)dr\Big|^{2}\Big]\\
 & \qquad\qquad=\sum_{k,l=1}^{n}\int_{\left(k-1\right)\triangle_{n}}^{k\triangle_{n}}\int_{\left(l-1\right)\triangle_{n}}^{l\triangle_{n}}\E\big[\left(f\left(X_{r}\right)-f\left(X_{(k-1)\triangle_{n}}\right)\right)\cdot\\
 & \qquad\qquad\qquad\qquad\qquad\qquad\qquad\qquad\cdot\left(f\left(X_{h}\right)-f\left(X_{(l-1)\triangle_{n}}\right)\right)\big]drdh.
\end{align*}
We bound the diagonal ($l=k)$ and off-diagonal ($l\neq k$) terms
separately. Consider first the diagonal case and $(k-1)\triangle_{n}\leq r\leq h\leq k\triangle_{n}$.
By the Markov property and stationarity of $X$ we can calculate the
expectation above explicitly. Indeed, we see that
\begin{align*}
 & \E\left[\left(f\left(X_{r}\right)-f\left(X_{(k-1)\triangle_{n}}\right)\right)\left(f\left(X_{h}\right)-f\left(X_{(k-1)\triangle_{n}}\right)\right)\right]\\
 & \qquad\qquad=\left\langle P_{h-r}f,f\right\rangle _{\mu}-\left\langle P_{h-\left(k-1\right)\triangle_{n}}f,f\right\rangle _{\mu}-\left\langle P_{r-\left(k-1\right)\triangle_{n}}f,f\right\rangle _{\mu}+\left\langle f,f\right\rangle _{\mu}\\
 & \qquad\qquad=\left\langle \left(P_{h-r}-I\right)f+\left(I-P_{h-\left(k-1\right)\triangle_{n}}\right)f+\left(I-P_{r-\left(k-1\right)\triangle_{n}}\right)f,f\right\rangle _{\mu}.
\end{align*}
Consequently, by symmetry in $r,h$
\begin{align*}
 & \sum_{k=1}^{n}\int_{\left(k-1\right)\triangle_{n}}^{k\triangle_{n}}\int_{\left(k-1\right)\triangle_{n}}^{k\triangle_{n}}\E\big[\left(f\left(X_{r}\right)-f\left(X_{(k-1)\triangle_{n}}\right)\right)\cdot\\
 & \qquad\qquad\qquad\qquad\qquad\qquad\cdot\left(f\left(X_{h}\right)-f\left(X_{(k-1)\triangle_{n}}\right)\right)\big]drdh\\
 & \qquad\qquad=2\sum_{k=1}^{n}\bigg<\bigg(\int_{\left(k-1\right)\triangle_{n}}^{k\triangle_{n}}\int_{\left(k-1\right)\triangle_{n}}^{h}\left(P_{h-r}-I\right)drdh\\
 & \qquad\qquad\qquad\qquad\qquad\qquad+\triangle_{n}\int_{\left(k-1\right)\triangle_{n}}^{k\triangle_{n}}\left(I-P_{h-\left(k-1\right)\triangle_{n}}\right)dh\bigg)f,f\bigg>_{\mu}\\
 & \qquad\qquad=2n\Big\langle\bigg(\int_{0}^{\triangle_{n}}\int_{0}^{h}\left(P_{h-r}-I\right)drdh+\triangle_{n}\int_{0}^{\triangle_{n}}\left(I-P_{h}\right)dh\bigg)f,f\Big\rangle_{\mu}.
\end{align*}
Since the generator $L$ is normal, by the functional calculus (see
\ref{sec:Appendix}) we can write this as 
\[
\left\langle \Psi(L)f,f\right\rangle _{\mu}=\int_{\sigma\left(L\right)}\Psi\left(\lambda\right)d\left\langle E_{\lambda}f,f\right\rangle _{\mu}
\]
for the measurable function
\[
\Psi(\lambda)=2n\left(\int_{0}^{\triangle_{n}}\int_{0}^{h}\left(e^{\lambda(h-r)}-1\right)drdh+\triangle_{n}\int_{0}^{\triangle_{n}}\left(1-e^{\lambda h}\right)dh\right),\,\,\,\,\lambda\in\mathbb{C}.
\]
Fix now $0\leq s\leq1$ such that for $z\in\{\lambda\in\C\,:\,Re(\lambda)\leq0\}$
we have $|1-e^{z}|\leq2|z|^{s}$. Since $L$ is the generator of a
Feller semigroup, we know that $\sigma(L)\subset\{\lambda\in\C\,:\,Re(\lambda)\leq0\}$.
We therefore conclude that $|\Psi(\lambda)|\leq4n\triangle_{n}^{2+s}|\lambda|^{s}$,
$\lambda\in\sigma(L)$. Hence the diagonal terms are bounded by
\begin{equation}
\int_{\sigma(L)}\left|\Psi(\lambda)\right|d\langle E_{\lambda}f,f\rangle_{\mu}\leq4T\triangle_{n}^{1+s}\int_{\sigma(L)}\left|\lambda\right|^{s}d\langle E_{\lambda}f,f\rangle_{\mu}=4T\triangle_{n}^{1+s}\|\left|L\right|{}^{s/2}f\|_{\mu}^{2},\label{eq:diagonalBound}
\end{equation}
which is true as long as $f\in\text{dom}(|L|^{s/2})$. For the off-diagonal
terms consider $(l-1)\triangle_{n}\leq r\leq(k-1)\triangle_{n}\leq h$.
Then, similar as before
\begin{align*}
 & \E\left[\left(f\left(X_{h}\right)-f\left(X_{\left(k-1\right)\triangle_{n}}\right)\right)\left(f\left(X_{r}\right)-f\left(X_{\left(l-1\right)\triangle_{n}}\right)\right)\right]\\
 & \qquad=\left\langle P_{h-r}f,f\right\rangle _{\mu}-\left\langle P_{h-\left(l-1\right)\triangle_{n}}f,f\right\rangle _{\mu}-\left\langle P_{\left(k-1\right)\triangle_{n}-r}f,f\right\rangle _{\mu}+\left\langle P_{\left(k-l\right)\triangle_{n}}f,f\right\rangle _{\mu}\\
 & \qquad=\left\langle P_{\left(k-1\right)\triangle_{n}-r}(P_{h-\left(k-1\right)\triangle_{n}}-I)(I-P_{r-\left(l-1\right)\triangle_{n}})f,f\right\rangle _{\mu}.
\end{align*}
If the generator is self-adjoint, then the $P_{u}$, $u\geq0$, are
as well. Then $P_{u}$ is positive, $P_{u}-I$ negative and $I-P_{u}$
again positive semidefinite. We can thus conclude that $P_{\left(k-1\right)\triangle_{n}-r}(P_{h-\left(k-1\right)\triangle_{n}}-I)(I-P_{r-\left(l-1\right)\triangle_{n}})$
is negative semidefinite and that the off-diagonal terms do not contribute
to the estimation error. In the more general case of a normal operator
$L$, we have to explicitly bound the off-diagonal terms. First, the
off-diagonal terms are equal to
\begin{align*}
 & 2\sum_{k>l=1}^{n}\int_{\left(k-1\right)\triangle_{n}}^{k\triangle_{n}}\int_{\left(l-1\right)\triangle_{n}}^{l\triangle_{n}}\E\big[\left(f\left(X_{r}\right)-f\left(X_{(k-1)\triangle_{n}}\right)\right)\cdot\\
 & \qquad\qquad\qquad\qquad\qquad\qquad\cdot\left(f\left(X_{h}\right)-f\left(X_{(l-1)\triangle_{n}}\right)\right)\big]drdh\\
 & \qquad=2\sum_{k>l=1}^{n}\bigg<\bigg(\int_{\left(k-1\right)\triangle_{n}}^{k\triangle_{n}}\int_{\left(l-1\right)\triangle_{n}}^{l\triangle_{n}}P_{\left(k-l\right)\triangle_{n}-\left(r-\left(l-1\right)\triangle_{n}\right)}\\
 & \qquad\qquad\qquad\qquad\Big(P_{h-\left(k-1\right)\triangle_{n}}-I\Big)\left(I-P_{r-\left(l-1\right)\triangle_{n}}\right)drdh\bigg)f,f\bigg>_{\mu}\\
 & \qquad=\bigg<2\bigg(\int_{0}^{\triangle_{n}}\int_{0}^{\triangle_{n}}\left(\sum_{k>l=1}^{n}P_{\left(k-l\right)\triangle_{n}-r}\right)\Big(P_{h}-I\Big)\left(I-P_{r}\right)drdh\bigg)f,f\bigg>_{\mu}\\
 & \qquad=\int_{\sigma\left(L\right)}\tilde{\Psi}\left(\lambda\right)d\left\langle E_{\lambda}f,f\right\rangle _{\mu},
\end{align*}
where we argue as before, but this time the function $\tilde{\Psi}$
is given by 
\[
\tilde{\Psi}(\lambda)=2\int_{0}^{\triangle_{n}}\int_{0}^{\triangle_{n}}\left(\sum_{k>l=1}^{n}e^{\lambda\left(\left(k-l\right)\triangle_{n}-r\right)\big)}\right)\Big(e^{\lambda h}-1\Big)\Big(1-e^{\lambda r}\Big)drdh,\,\,\,\,\lambda\in\mathbb{C}.
\]
We will show that there exists a universal constant $\tilde{C}<\infty$
such that 
\begin{equation}
\left|\tilde{\Psi}\left(\lambda\right)\right|\leq\tilde{C}T|\lambda|^{s}\triangle_{n}^{1+s},\,\,\,\,\lambda\in\sigma(L).\label{eq:internal}
\end{equation}
As in \eqref{eq:diagonalBound} we then can bound the off-diagonal
terms by 
\[
\int_{\sigma(L)}\left|\Psi(\lambda)\right|d\langle E_{\lambda}f,f\rangle_{\mu}\leq\tilde{C}T\triangle_{n}^{1+s}\int_{\sigma(L)}\left|\lambda\right|^{s}d\langle E_{\lambda}f,f\rangle_{\mu}=\tilde{C}T\|\left|L\right|{}^{s/2}f\|_{\mu}^{2}\triangle_{n}^{1+s}
\]
for $f\in\text{dom}(|L|^{s/2})$. Combining this with \eqref{eq:diagonalBound}
yields the claim. For $\lambda=0$ we also have $\tilde{\Psi}(\lambda)=0$.
It is therefore sufficient to consider $\lambda\neq0$. In order to
bound $\tilde{\Psi}$ in that case we calculate
\[
\sum_{k>l=1}^{n}e^{\lambda\left(k-l-1\right)\triangle_{n}}=\sum_{l=1}^{n}\frac{1-e^{\lambda\left(n-l\right)\triangle_{n}}}{1-e^{\lambda\triangle_{n}}}=\frac{n}{1-e^{\lambda\triangle_{n}}}-\frac{1-e^{\lambda n\triangle_{n}}}{\left(1-e^{\lambda\triangle_{n}}\right)^{2}}.
\]
Hence
\begin{align*}
\left|\triangle_{n}^{2}\left(1-e^{\lambda\triangle_{n}}\right)^{2}\sum_{k>l=1}^{n}e^{\lambda\left(k-l-1\right)\triangle_{n}\big)}\right| & \leq2n\triangle_{n}^{2}\left|\lambda\triangle_{n}\right|^{s}+2\triangle_{n}^{2}\left|\lambda n\triangle_{n}\right|^{s}\\
 & \leq4T\triangle_{n}^{1+s}\left|\lambda\right|^{s}.
\end{align*}
Therefore \eqref{eq:internal} follows if we can show that
\begin{equation}
\left|\triangle_{n}^{-2}\left(1-e^{\lambda\triangle_{n}}\right)^{-2}\int_{0}^{\triangle_{n}}\int_{0}^{\triangle_{n}}e^{\lambda\left(\triangle_{n}-r\right)}\Big(e^{\lambda h}-1\Big)\Big(1-e^{\lambda r}\Big)drdh\right|\label{eq:thirdBound}
\end{equation}
is bounded by a universal constant. To show this, let $z=\lambda\triangle_{n}$
and note that
\begin{equation}
\triangle_{n}^{-1}\left|\frac{\int_{0}^{\triangle_{n}}\left(1-e^{\lambda h}\right)dh}{1-e^{\lambda\triangle_{n}}}\right|=\triangle_{n}^{-1}\left|\frac{\triangle_{n}}{1-e^{\lambda\triangle_{n}}}-\lambda^{-1}\right|=\left|\frac{1}{z}-\frac{1}{e^{z}-1}\right|,\label{eq:firstBound}
\end{equation}
\begin{equation}
\triangle_{n}^{-1}\left|\int_{0}^{\triangle_{n}}\frac{\Big(e^{\lambda\left(\triangle_{n}-r\right)}-e^{\lambda\triangle_{n}}\Big)}{1-e^{\lambda\triangle_{n}}}dr\right|=\left|-\frac{1}{z}-\frac{e^{z}}{1-e^{z}}\right|=\left|\frac{1}{z}-\frac{1}{e^{z}-1}+1\right|.\label{eq:secondBound}
\end{equation}
\eqref{eq:firstBound} converges to $1$ and \eqref{eq:secondBound}
converges to $0$ as $|z|\rightarrow\infty$. We further see that
\eqref{eq:firstBound} converges to $1/2$ and \eqref{eq:secondBound}
converges to $3/2$ when $|z|\rightarrow0$ and when considering only
$z\in\{\lambda\in\mathbb{C}:\text{Re}(\lambda)\leq0\}$. This implies
a universal constant bounding \eqref{eq:thirdBound}. This proves
\eqref{eq:internal} and thus the claim.
\end{proof}

\subsection{Proof of Corollary \ref{cor:Not stationary initial distribution}}
\begin{proof}
Recall that $X_{0}\overset{d}{\sim}\eta$ and $\eta$ is absolutely
continuous with respect to $\mu$ with bounded density $d\eta/d\mu$.
We indicate by $\E_{\nu}$ for a measure $\nu$ the initial distribution
of $X$. By conditioning on $X_{0}$ we see that
\begin{align*}
\left\Vert \Gamma_{T}(f)-\hat{\Gamma}_{T,n}(f)\right\Vert _{L^{2}(\P)}^{2} & =\E_{\eta}\Big[\Big|\Gamma_{T}(f)-\hat{\Gamma}_{T,n}(f)\Big|^{2}\big|\Big]\\
 & =\int_{\mathcal{S}}\E\Big[\Big|\Gamma_{T}(f)-\hat{\Gamma}_{T,n}(f)\Big|^{2}\big|X_{0}\Big]_{X_{0}=x}d\eta(x)\\
 & \leq\Big\|\frac{d\eta}{d\mu}\Big\|_{\infty,\mu}\int_{\mathcal{S}}\E\Big[\Big|\Gamma_{T}(f)-\hat{\Gamma}_{T,n}(f)\Big|^{2}\big|X_{0}\Big]_{X_{0}=x}d\mu(x)\\
 & =\Big\|\frac{d\eta}{d\mu}\Big\|_{\infty,\mu}\E_{\mu}\Big[\Big|\Gamma_{T}(f)-\hat{\Gamma}_{T,n}(f)\Big|^{2}\big|\Big].
\end{align*}
The conclusion follows immediately from Theorem \ref{thm:OT:Rates for Bessel potential spaces}.
\end{proof}

\subsection{Proof of Theorem \ref{thm:Convergence of the occupation measure to the stationary measure}}
\begin{proof}
By the triangle inequality we have for $f\in L^{2}(\mu)$ and Theorem
\ref{thm:OT:Rates for Bessel potential spaces} that
\begin{align*}
 & \left\Vert T^{-1}\hat{\Gamma}_{T,n}\left(f\right)-\int_{\mathcal{S}}f(x)\,d\mu(x)\right\Vert _{L^{2}(\P)}\\
 & \qquad\leq T^{-1}\left\Vert \hat{\Gamma}_{T,n}\left(f\right)-\Gamma_{T}\left(f\right)\right\Vert _{L^{2}(\P)}+\left\Vert T^{-1}\Gamma_{T}\left(f\right)-\int_{\mathcal{S}}f(x)\,d\mu(x)\right\Vert _{L^{2}(\P)}\\
 & \qquad\leq\frac{C}{\sqrt{T}}\|f\|_{L^{2}(\mu)}\triangle_{n}^{1/2}+\left\Vert T^{-1}\Gamma_{T}\left(f\right)-\int_{\mathcal{S}}f(x)\,d\mu(x)\right\Vert _{L^{2}(\P)}
\end{align*}
for a universal constant $C$. The claimed bound for the second term
is well-known, but we give the proof here to complement the proof
of Theorem \ref{thm:OT:Rates for Bessel potential spaces}. Consider
$f$ such that $f_{0}=f-\int_{\mathcal{}}fd\mu\in\text{dom}(|L|^{-1/2})$.
By linearity of the occupation time functional it holds
\[
T^{-1}\Gamma_{T}(f)-\int_{\mathcal{S}}fd\mu=T^{-1}\Gamma_{T}(f_{0}).
\]
Fubini's theorem yields
\begin{align*}
\E\Big[\Big|T^{-1}\Gamma_{T}\left(f_{0}\right)\Big|^{2}\Big] & =T^{-2}\int_{0}^{T}\int_{0}^{T}\E\left[f_{0}\left(X_{r}\right)f_{0}\left(X_{h}\right)\right]drdh\\
 & =2T^{-2}\int_{0}^{T}\int_{0}^{h}\left\langle P_{h-r}f_{0},f_{0}\right\rangle _{\mu}drdh\\
 & =\int_{\sigma(L)}\Psi(\lambda)\,d\langle E_{\lambda}f_{0},f_{0}\rangle_{\mu},
\end{align*}
where the measurable function $\Psi$ is defined by
\begin{align*}
\Psi(\lambda) & =2T^{-2}\int_{0}^{T}\int_{0}^{h}e^{\lambda(h-r)}drdh=2\frac{e^{\lambda T}-1-\lambda T}{\lambda^{2}T^{2}}=2\frac{(\lambda T)^{-1}(e^{\lambda T}-1)-1}{\lambda T},
\end{align*}
and where $\Psi(0)=1$ by continuous extension. Since $z\to z^{-1}(e^{z}-1)-1$
is bounded on the left half-plane $\{z\in\C\,:\,\text{Re}(z)\leq0\}$,
there exists a constant $\tilde{C}<\infty$ such that
\[
\left|\Psi\left(\lambda\right)\right|\leq\frac{\tilde{C}}{\left|\lambda\right|T},\quad\lambda\in\sigma(L).
\]
We conclude that
\[
\E\Big[\Big|T^{-1}\Gamma_{T}\left(f_{0}\right)\Big|^{2}\Big]\leq\frac{C}{T}\int_{\sigma(L)}\left|\lambda\right|^{-1}d\left\langle E_{\lambda}f_{0},f_{0}\right\rangle _{\mu}=\frac{C}{T}\|\left|L\right|^{-1/2}f_{0}\|_{\mu}^{2}.\tag*{{\qedhere}}
\]
\end{proof}

\subsection{Proof of Theorem \ref{thm:BrownianMotion}}
\begin{proof}
We consider only the case that $f\in H^{s}(\R^{d})$. The case of
an indicator function is analogous. The main idea of the proof is
to approximate the process $X$ by reflected processes for which Theorem
\ref{thm:OT:Rate for Sobolev functions} can be applied. Choose $M$
big enough such that $\operatorname{supp}(\eta)\subset[-M,M]$ and
let $\tau_{M}=\inf\{r>0\,:\,|X_{r}|>M\}$ the first time $X$ exits
from $[-M,M]$. By dominated convergence, for any fixed $T>0$, it
holds that
\begin{equation}
\left\Vert \Gamma_{T}(f)-\hat{\Gamma}_{T,n}(f)\right\Vert _{L^{2}(\P)}=\lim_{M\to\infty}\left\Vert \big(\Gamma_{T}(f)-\hat{\Gamma}_{T,n}(f)\big)\I(T<\tau_{M})\right\Vert _{L^{2}(\P)}.\label{eq:OT:approximation by reflected processes}
\end{equation}
We will show that there exists a universal constant $C$, independent
of $M$, such that 
\begin{equation}
\left\Vert \big(\Gamma_{T}(f)-\hat{\Gamma}_{T,n}(f)\big)\I(t<\tau_{M})\right\Vert _{L^{2}(\P)}\leq C\Big\|\frac{d\eta}{dx}\Big\|_{\infty,\mu}^{1/2}\|f\|_{H^{s}}\sqrt{T}\triangle_{n}^{\frac{1+s}{2}}.\label{eq:Bound uniform in M}
\end{equation}
For this define
\[
f_{M}(x)=\begin{cases}
x-4kM & :\,(4k-1)M\leq x<(4k+1)M\\
(4k+2)M-x & :\,(4k+1)M\leq x<(4k+3)M
\end{cases}.
\]
Applying the Itô-Tanaka formula one can check that the process $X_{r}^{(M)}:=f_{M}(X_{r})$
is a reflected Brownian motion with barriers at $-M,M,$ i.e. $(X_{r}^{(M)})_{r\geq0}$
satisfies \eqref{eq:reflectedDiffusion} with $L=-M$, $R=M$, $b=0$,
$\sigma=1$ and some process $(K_{r})_{r\geq0}$. See \citet[Chapter I.23]{gihman2015stochastic}
for a similar construction of reflected diffusion processes.  Furthermore,
for all $0\leq r\leq T\leq\tau_{M}$ we have
\begin{equation}
B_{r}=X_{r}^{(M)}.\label{eq:OT:internal 2}
\end{equation}
In the following denote by $\Gamma_{T}^{(M)},\hat{\Gamma}_{T,n}^{(M)}$
the integral functional and the Riemann estimator with respect to
the process $(X_{r}^{(M)})_{r\geq0}$. From \eqref{eq:OT:internal 2}
we obtain that
\begin{align*}
\left\Vert \big|\Gamma_{T}(f)-\hat{\Gamma}_{T,n}(f)\big|\I(T<\tau_{M})\right\Vert _{L^{2}(\P)} & =\left\Vert \big|\Gamma_{T}^{(M)}-\hat{\Gamma}_{T,n}^{(M)}(f)\big|\I(T<\tau_{M})\right\Vert _{L^{2}(\P)}\\
 & \leq\left\Vert \Gamma_{T}^{(M)}-\hat{\Gamma}_{T,n}^{(M)}(f)\right\Vert _{L^{2}(\P)}.
\end{align*}
Since $\operatorname{supp}(\eta)\subset[-M,M]$, in particular $X_{0}=X_{0}^{(M)}$,
i.e. $X_{0}^{(M)}$ has distribution $\eta.$ Note that $\eta$ is
in general not the stationary distribution of $(X_{r}^{(M)})_{r\geq0}$.
From Example \eqref{ex:OT:Scalar-diffusion-process with reflection}
we actually know that the stationary distribution $\mu_{M}$ of the
reflected Brownian motion on $[-M,M]$ has Lebesgue density $d\mu_{M}/dx=(2M)^{-1}\I_{[-M,M]}$.
Theorem \ref{thm:OT:Rate for Sobolev functions}, with $A=1$, therefore
implies together with the argument from the proof of Corollary \ref{cor:Not stationary initial distribution}
that 
\begin{align*}
\left\Vert \Gamma_{T}^{(M)}-\hat{\Gamma}_{T,n}^{(M)}(f)\right\Vert _{L^{2}(\P)} & \leq C\Big\|\frac{d\mu_{M}}{dx}\Big\|_{\infty,\mu}^{1/2}\Big\|\frac{d\eta}{d\mu_{M}}\Big\|_{\infty,\eta}^{1/2}\|f\|_{H^{s}}\sqrt{T}\triangle_{n}^{\frac{1+s}{2}}.
\end{align*}
Observe that 
\[
\int gd\eta=\int_{-M}^{M}g\frac{d\eta}{dx}dx=2M\int g\frac{d\eta}{dx}d\mu_{M}
\]
for any bounded continuous function $g$, i.e. $d\eta/d\mu_{M}=2M\,d\eta/dx$.
We conclude that 
\begin{align*}
\left\Vert \Gamma_{T}^{(M)}-\hat{\Gamma}_{T,n}^{(M)}(f)\right\Vert _{L^{2}(\P)} & \leq C\Big\|\frac{d\eta}{dx}\Big\|_{\infty,\eta}^{1/2}\|f\|_{H^{s}}\sqrt{T}\triangle_{n}^{\frac{1+s}{2}},
\end{align*}
which is \eqref{eq:Bound uniform in M}.
\end{proof}

\appendix

\section*{\label{sec:Appendix}Appendix. A short review of semigroup theory
and functional calculus}

We will briefly recall the basic objects needed in the theory of semigroups
and the functional calculus for normal operators. For more details
see \citet{engel1999one} and \citet{rudin2006functional}. Let $\mu$
be any probability measure on $(\c S,\c B(\c S))$. On the induced
Hilbert space $L^{2}(\mu)$, denote by $(P_{r})_{r\geq0}$ the Markov
semigroup associated with $X$ which satisfies $P_{r}f(x)=\E[f(X_{r})|X_{0}=x]$
for $f\in L^{2}(\mu)$, $x\in\c S$ and $P_{r+s}=P_{r}P_{s}$, $r,s\geq0$.
The infinitesimal generator of the semigroup is defined as 
\[
Lf=\lim_{r\rightarrow0}\frac{P_{r}f-f}{r},\,\,\,\,f\in\text{dom}(L),
\]
where $\text{dom}(L)\subset L^{2}(\mu)$ is the set of all functions
$f$ for which this limit exists. If $(P_{r})_{r\geq0}$ is strongly
continuous, i.e. $P_{r}f\r{}f$ in $L^{2}(\mu)$ as $r\rightarrow0$
for all $f\in L^{2}(\mu)$, then the semigroup is Feller. This is
true for most Markov processes in practice, including Lévy processes
and many diffusions. In the Feller case, $L$ is a densely defined
closed linear and usually unbounded operator on its domain with spectrum
$\sigma(L)\subset\{\lambda\in\mathbb{C}:\text{Re}(\lambda)\leq0\}$.
In order to define fractional powers of the generator we further need
the operator $L$ (and thus the operators $P_{r}$) to be normal,
i.e. $LL^{*}=L^{*}L$ where $L^{*}$ is the Hilbert space adjoint
of $L$. In that case the spectral theorem (Theorem 13.33 of \citet{rudin2006functional})
guarantees the existence of a resolution of the identity or spectral
measure $(E_{A})_{A\in\c B(\mathbb{C})}$ on $L^{2}(\mu)$. This means
that $(E_{A})_{A\in\c B(\mathbb{C})}$ is a family of orthogonal projections
$E_{A}:L^{2}(\mu)\to L^{2}(\mu)$ for Borel sets $A\subset\mathbb{C}$
such that for every $f,g\in L^{2}(\mu)$ the map $A\mapsto\left\langle E_{A}f,g\right\rangle _{\mu}$
is a complex measure, supported on $\sigma(L)$. Moreover, $A\mapsto\left\langle E_{A}f,f\right\rangle _{\mu}$
is a positive measure with total variation $\left\langle E_{\mathbb{C}}f,f\right\rangle _{\mu}=\norm f_{\mu}$.
By the spectral theorem we can associate to any measurable function
$\Psi:\mathbb{C}\mapsto\mathbb{C}$ a densely defined closed operator
$\Psi(L)$ with domain $\text{dom}(\Psi(L)):=\{f\in L^{2}(\mu):\int_{\sigma(L)}|\Psi(\lambda)|^{2}d\langle E_{\lambda}f,f\rangle_{\mu}<\infty\}$
by the relation
\[
\left\langle \Psi\left(L\right)f,g\right\rangle _{\mu}=\int_{\sigma\left(L\right)}\Psi\left(\lambda\right)d\langle E_{\lambda}f,g\rangle_{\mu},\,\,\,\,f\in\text{dom}(\Psi\left(L\right)),g\in L^{2}\left(\mu\right).
\]
It satisfies $\norm{\Psi(L)f}_{\mu}=\int_{\sigma(L)}|\Psi(\lambda)|^{2}d\left\langle E_{\lambda}f,f\right\rangle _{\mu}$.
In particular, we can define the fractional operators $|L|^{s/2}$
on $\c D^{s}(L):=\text{dom}(|L|{}^{s/2})$, for $0\leq s\leq1$. By
the spectral theorem for normal semigroups (Theorem 13.38 of \citet{rudin2006functional})
we realize the semigroup in its usual exponential form, i.e. $P_{r}=\Psi(L)$
with $\Psi(x)=e^{rx}$, $r\geq0$. 

\section*{Acknowledgements}

Both authors gratefully acknowledge the financial support of the DFG
Research Training Group 1845 ”Stochastic Analysis with Applications
in Biology, Finance and Physics”.

\bibliographystyle{apalike}
\bibliography{bibliography}

\begin{thebibliography}{}

\bibitem[Adams and Fournier, 2003]{adams2003sobolev}
Adams, R.~A. and Fournier, J. J.~F. (2003).
\newblock {\em {Sobolev Spaces}}.
\newblock Pure and Applied Mathematics. Elsevier Science.

\bibitem[Bakry et~al., 2013]{bakry2013analysis}
Bakry, D., Gentil, I., and Ledoux, M. (2013).
\newblock {\em {Analysis and Geometry of Markov Diffusion Operators}}.
\newblock Grundlehren der mathematischen Wissenschaften. Springer International
  Publishing.

\bibitem[Bass, 2006]{bass2006diffusions}
Bass, R.~F. (2006).
\newblock {\em {Diffusions and Elliptic Operators}}.
\newblock Probability and Its Applications. Springer New York.

\bibitem[Catellier and Gubinelli, 2016]{Catellier2016}
Catellier, R. and Gubinelli, M. (2016).
\newblock {Averaging along irregular curves and regularisation of ODEs}.
\newblock {\em Stochastic Processes and their Applications}, 126(8):2323--2366.

\bibitem[Chen, 2006]{chen2006eigenvalues}
Chen, M.~F. (2006).
\newblock {\em {Eigenvalues, Inequalities, and Ergodic Theory}}.
\newblock Probability and Its Applications. Springer London.

\bibitem[Chesney et~al., 1997]{Chesney1997}
Chesney, M., Jeanblanc-Picqu{\'{e}}, M., and Yor, M. (1997).
\newblock {Brownian Excursions and Parisian Barrier Options}.
\newblock {\em Source: Advances in Applied Probability Adv. Appl. Prob},
  29(29).

\bibitem[Chojnowska-Michalik and Goldys, 2002]{Chojnowska-Michalik2002}
Chojnowska-Michalik, A. and Goldys, B. (2002).
\newblock {Symmetric Ornstein-Uhlenbeck semigroups and their generators}.
\newblock {\em Probability Theory and Related Fields}, 124(4):459--486.

\bibitem[Chorowski, 2015]{Chorowski2015a}
Chorowski, J. (2015).
\newblock {Nonparametric volatility estimation in scalar diffusions: Optimality
  across observation frequencies}.
\newblock {\em arXiv:1507.07139}.

\bibitem[Diehl et~al., 2016]{Diehl2016a}
Diehl, J., Gubinelli, M., and Perkowski, N. (2016).
\newblock {The Kardar-Parisi-Zhang equation as scaling limit of weakly
  asymmetric interacting Brownian motions}.
\newblock {\em arXiv:1606.02331}.

\bibitem[Dion and Genon-Catalot, 2016]{Dion2016}
Dion, C. and Genon-Catalot, V. (2016).
\newblock {Bidimensional random effect estimation in mixed stochastic
  differential model}.
\newblock {\em Statistical Inference for Stochastic Processes}, 19(2):131--158.

\bibitem[Engel and Nagel, 1999]{engel1999one}
Engel, K.-J. and Nagel, R. (1999).
\newblock {\em {One-Parameter Semigroups for Linear Evolution Equations}}.
\newblock Graduate Texts in Mathematics. Springer New York.

\bibitem[Ethier and Kurtz, 1986]{Ethier1986}
Ethier, S.~N. and Kurtz, T.~G. (1986).
\newblock {\em {Markov processes : characterization and convergence}}.
\newblock Wiley.

\bibitem[Ganychenko, 2015]{Ganychenko2015a}
Ganychenko, I. (2015).
\newblock {Fast L 2 -approximation of integral-type functionals of Markov
  processes}.
\newblock {\em Modern Stochastics: Theory and Applications}, 2:165--171.

\bibitem[Ganychenko et~al., 2015]{Ganychenko2015}
Ganychenko, I.~V., Knopova, V.~P., and Kulik, A.~M. (2015).
\newblock {Accuracy of discrete approximation for integral functionals of
  Markov processes}.
\newblock {\em Modern Stochastics: Theory and Applications}, 2(4):401--420.

\bibitem[Gihman and Skorohod, 2015]{gihman2015stochastic}
Gihman, I. and Skorohod, A. (2015).
\newblock {\em {Stochastic Differential Equations}}.
\newblock Springer.

\bibitem[Gobet and Matulewicz, 2016]{Gobet2016}
Gobet, E. and Matulewicz, G. (2016).
\newblock {Parameter estimation of Ornstein--Uhlenbeck process generating a
  stochastic graph}.
\newblock {\em Statistical Inference for Stochastic Processes}, pages 1--25.

\bibitem[Hansen et~al., 1998]{Hansen1998}
Hansen, L.~P., {Alexandre Scheinkman}, J., and Touzi, N. (1998).
\newblock {Spectral methods for identifying scalar diffusions}.
\newblock {\em Journal of Econometrics}, 86(1):1--32.

\bibitem[Hoffmann, 1999]{Hoffmann1999}
Hoffmann, M. (1999).
\newblock {L p Estimation of the Diffusion Coefficient}.
\newblock {\em Bernoulli}, 5(3):447.

\bibitem[Hugonnier, 1999]{hugonnier1999feynman}
Hugonnier, J.-N. (1999).
\newblock {The Feynman--Kac formula and pricing occupation time derivatives}.
\newblock {\em International Journal of Theoretical and Applied Finance},
  2(02):153--178.

\bibitem[Kipnis and Varadhan, 1986]{Kipnis1986}
Kipnis, C. and Varadhan, S. R.~S. (1986).
\newblock {Central limit theorem for additive functionals of reversible Markov
  processes and applications to simple exclusions}.
\newblock {\em Communications in Mathematical Physics}, 104(1):1--19.

\bibitem[Kohatsu-Higa et~al., 2014]{Kohatsu-Higa2014}
Kohatsu-Higa, A., Makhlouf, A., and Ngo, H. (2014).
\newblock {Approximations of non-smooth integral type functionals of one
  dimensional diffusion processes}.
\newblock {\em Stochastic Processes and their Applications}, 124(5):1881--1909.

\bibitem[Ngo and Ogawa, 2011]{Ngo2011}
Ngo, H.-L. and Ogawa, S. (2011).
\newblock {On the discrete approximation of occupation time of diffusion
  processes}.
\newblock {\em Electronic Journal of Statistics}, 5:1374--1393.

\bibitem[Pavliotis, 2014]{pavliotis2014stochastic}
Pavliotis, G.~A. (2014).
\newblock {\em {Stochastic Processes and Applications: Diffusion Processes, the
  Fokker-Planck and Langevin Equations}}.
\newblock Texts in Applied Mathematics. Springer New York.

\bibitem[Pollett, 2003]{Pollett2003}
Pollett, P. (2003).
\newblock {Integrals for continuous-time Markov chains}.
\newblock {\em Mathematical Biosciences}, 182(2):213--225.

\bibitem[Rudin, 2006]{rudin2006functional}
Rudin, W. (2006).
\newblock {\em {Functional Analysis}}.
\newblock International series in pure and applied mathematics. McGraw-Hill.

\bibitem[{V. Gol'dshtein}, 2009]{10.2307/40302922}
{V. Gol'dshtein}, A.~U. (2009).
\newblock {Weighted Sobolev Spaces and Embedding Theorems}.
\newblock {\em Transactions of the American Mathematical Society},
  361(7):3829--3850.

\bibitem[Watanabe, 1984]{watanabe1984lectures}
Watanabe, S. (1984).
\newblock {\em {Lectures on Stochastic Differential Equations and Malliavin
  Calculus}}.
\newblock Lectures on mathematics and physics / Tata institute of fundamental
  research: Mathematics. Springer.

\end{thebibliography}

\end{document}